\documentclass[12pt]{elsarticle}

\usepackage{graphicx}
\usepackage{amssymb}

\usepackage{amsmath}
\usepackage{amsfonts}
\usepackage{graphicx}
\usepackage{bm}
\usepackage{relsize}
\usepackage{float} 
\usepackage{caption,setspace}
\usepackage{booktabs}
\usepackage{multirow}
\usepackage{subcaption}
\usepackage{textcomp}
\usepackage{cleveref}
\usepackage{breqn}
\usepackage{tikz}
\usetikzlibrary{math}
\usetikzlibrary{shapes.geometric, decorations.pathreplacing}

\makeatletter
\def\ps@pprintTitle{%
  \let\@oddhead\@empty
  \let\@evenhead\@empty
  \let\@oddfoot\@empty
  \let\@evenfoot\@oddfoot
}
\makeatother

\begin{document}


\begin{frontmatter}


\title{Curvilinear High-Order Mimetic Differences that Satisfy Conservation Laws}



\author[1]{Angel Boada}
\ead{aboadavelazco@sdsu.edu} 
\author[1,2]{Johnny Corbino}
\ead{jcorbino@lbl.gov}
\author[1]{Miguel Dumett}
\ead{mdumett@sdsu.edu}
\author[1]{Jose Castillo}
\ead{jcastillo@sdsu.edu}
\address[1]{Computational Science Research Center, 5500 Campanile Drive, San Diego, CA 92182}
\address[2]{Lawrence Berkeley National Laboratory, Berkeley, CA 94720}

\begin{abstract}

We investigate the construction and usage of mimetic operators in curvilinear staggered grids. Specifically, we extend the Corbino-Castillo operators so they can be utilized to solve problems in non-trivial geometries. We prove that the resulting curvilinear operators satisfy the discrete analog of the extended Gauss-Divergence theorem. In addition, we demonstrate energy and mass conservation in curvilinear coordinates for the acoustic wave equation. These findings are illustrated in two-dimensional and three-dimensional elliptic/hyperbolic equations and can be extended to other partial differential equations as well. 
\end{abstract}


\end{frontmatter}


\section{Introduction}
\label{S:1}


Mimetic methods are discrete numerical schemes based on mimetic difference operators. These operators are discrete analogs of spatial differential vector calculus operators like divergence, gradient, curl, Laplacian, etc.  They satisfy, in the discrete sense, the same properties that the continuum operators do. They mimic solution symmetries, conservation laws, vector calculus identities, and other important properties of continuum partial differential equations (PDE) mathematical models. These mimetic operators are constructed on staggered grids. 

Early work in ``mimetic methods'' was developed by A.A. Samarskii at Moscow State University during the middle of the twentieth century \cite{veiga, lipnikov, rodrigo}. Because of its implementation simplicity, mimetic methods were built initially utilizing finite differences.  These techniques were originally translated into English as support-operators methods \cite{shashkov}. These approaches have been utilized in a plethora of applications since then by a team of researchers mainly at Los Alamos National Laboratory (LANL) \cite{veiga, shashkov}. 

Support-operator techniques proceed in the following manner:
\begin{enumerate}
\item Write the PDE in terms of first-order invariant operators (divergence, gradient, and curl).
\item Select the properties to mimic. These integral laws are mainly Green identity type relationships based on line, surface, and volume integrals equivalence formulas (namely Green's, Gauss', and Stokes vector calculus theorems) of the divergence, gradient and curl operators. For example, in \cite{shashkov}, even though the book refers to five relationships, all mimetic discrete analogs for conservation laws are derived utilizing only the first one, specifically the extended Gauss divergence theorem (as in the high-order mimetic differences of Castillo and Grone \cite{Grone} and Corbino-Castillo \cite{Corbino}).
\item Select a subset of points of a non-staggered grid, where each scalar, vector, and tensor fields will be defined.
\item Choose one operator among the divergence, gradient or curl and define it (called the prime operator).
\item All other operators (referred as derived operators) are derived from the integral laws that relates the prime operator with each of the derived ones.
\item For two-dimensions (2D) and two-dimensions (3D), the divergence, gradient and curl operators are defined at nodes or cell centers.
\end{enumerate}
Even though, discrete analogs of the divergence and gradient operators are dual to each other in the presence of zero boundary conditions (which triggers symmetric Laplacian operators), the support-operator methods utilize no weight matrices for approximating integrals (inner products). It has been proven in \cite{Grone} and glimpsed at \cite{Kreiss} that high-order accurate discrete analogs of the first-order invariant differential operators are not available unless one utilizes positive definite weighted inner products. In this case, when scalar and/or vector fields are zero at the boundary, the duality relationship between gradient and divergence discrete analogs still holds for the weighted inner products \cite{Grone, Miranda, Corbino}.

Aware of this limitation of the support-operator method, in 2003 Castillo and Grone \cite{Grone} introduced high-order \emph{mimetic operators} (which we refer to high-order mimetic method, or mimetic differences for short, to distinguish from the support-operator method) that are uniformly accurate on staggered grids, with inner products given by positive diagonal weights and whose gradient, divergence and curl operators are given by parameterized families of operators, with the number of parameters varying accordingly to the degree of accuracy of the corresponding operators. Later, in 2019, Corbino-Castillo \cite{Corbino} found a unique (without parameters) similar high-order mimetic difference operators satisfying the same properties of the Castillo-Grone operators that are compact, and that have shorter band. In addition, MOLE (mimetic operator library enhanced), an open source library for the Corbino-Castillo operators, was developed at the time by J. Corbino \cite{altmolerepo}.

Such mimetic difference operators were constructed specifically to satisfy high-order discrete analogs of the fundamental theorem of calculus and the integration by parts formula \cite{Grone, Miranda}. Later, with the introduction of high-order mimetic interpolation operators \cite{Dumett-Castillo-Interpolators}, it was possible to shown that these operators satisfy the extended Gauss-Divergence theorem as well as fundamental vector calculus identities \cite{Dumett-Castillo-Identities}. These operators are built on a staggered grid. Moreover, these operators have a uniform order of accuracy in the boundary as well as the interior of the problem domain. The discrete mimetic difference operators, D and G, are constructed independently of each other as is the case of the divergence and gradient in vector calculus.

When compared to the support-operator methodology 1-6 above, mimetic differences write the PDE in terms of gradient, divergence, and curl (defined in terms of the divergence) operators, which are defined to act, respectively, on centers and faces (and not to act on nodes), independently of the PDE, and which by construction, like the support-operator methods, satisfy a discrete analog of the extended Gauss divergence theorem. Moreover, the 2D and 3D versions of these operators are conveniently generated as Kronecker products of the respective one-dimensional (1D) versions of the operator and some identity matrices \cite{Miranda,Dumett-Castillo-Identities}.

Nevertheless, during the last two decades, mimetic methods have been extended to utilize a framework similar to finite elements \cite{veiga}. This new method is called mimetic finite difference method (MFD) and it is the result of a collaboration between researchers at LANL and the Universit\`a degli Studi di Milano (Italy). Moreover, MFD has gained interest from several investigators of Europe and the US, who have recently found new connections between MFD and finite elements, especially for applications in electromagnetism \cite{adler, rodrigo}. 

Analogously as the support-operator method, MFD relies on coordinate invariant formulations based on different forms of the Stokes theorem for their primary mimetic gradient, divergence and curl operators (not one, as in the support-operator method, but each of them). With this set of primary first-order operators, a second set of gradient, divergence and curl operators (the duality operators) are created via duality relationships with proper inner products and the corresponding Green's identities of the Stokes' invariant formulations. However, the main difference of MFD with the support-operator method and mimetic differences is that instead of considering the integral identities for PDE classical solutions, MFD considers weak solutions \cite{evans}. Furthermore, MFD follows a finite element approach: it considers four different types of discrete fields defined by degrees of freedom associated to vertices, edges, faces, and elements. MFD is able to define all operators (including inner products) in an element-wise fashion: the invariant formulations are satisfied at each element and the whole domain mass and stiff matrices are built as a sum over all elements. 

This construction allows continuous first-order differential operators and their corresponding primal MFD versions to satisfy a de Rham cohomology diagram \cite{bredon}, from which kernels (that satisfy expected properties) of the primal and dual operators can be obtained (the former also satisfied by Castillo-Grone and Corbino-Castillo mimetic differences). In addition, second-order operators (Laplacian and curl-curl) are obtained by composing a dual and a primal operator. Furthermore, reconstruction operators are introduced as general inverses of projection operators (from the continuum domain to the discrete domain) to define consistency and stability conditions of this method, issues that are not addressed directly neither by the support-operator method nor by mimetic differences.

The choice of working on each element separately, in a finite element fashion, allows MFD to work on 1D, 2D and 3D different kind of grids, even on unstructured ones (however, in the latter case, some of the MFD identities do not hold and/or MFD operators are not able to be defined) \cite{veiga}. Moreover, high-order accuracy is achieved by increasing the number of degrees of freedom utilized in each of the four topologically different elements, in a finite element fashion, greatly increasing the number of unknowns that the discrete system has to solved for.

In our approach, mimetic difference operators are defined over computational grids that are Cartesian products of intervals (2D and 3D operators are defined via Kronecker versions of the 1D operators and convenient identity matrices). This, in principle, limits the domain complexity of the domain geometry that mimetic differences can solve \cite{delapuente}. One the goals of this paper is to extend the mimetic differences to curvilinear grids, increasing significantly the number of domains that mimetic differences can target. On the other hand, the mimetic differences high-order accuracy is not worked out element by element but by through several elements or cells, especially near the boundary grid elements. This way high-order accuracy is also guaranteed at any element. As a consequence of this fact, high-order mimetic difference operators have a relatively larger band than their corresponding mimetic difference lower accuracy counterparts and they have an explicit form that depends only on the number of cells of the domain.

In spite of the history of mimetic methods, one can find in the literature that the idea of creating numerical schemes that preserve some continuum vector calculus properties can be also found in other discretization methods, for example summation by parts \cite{ranocha}, finite volume \cite{leveque}, mixed finite element method \cite{boffi}, and an exposition of that is out of the scope of this paper. 

The present document falls under the approach of Corbino-Castillo mimetic differences but it can easily be adapted to the Castillo-Grone approach.

Recent applications of these mimetic operators involve highly non-linear problems like Richardson's equation for unsaturated fluid flow \cite{velazco2020high}, and elliptic problems in heterogeneous anisotropic media \cite{boada2020high}. However, to accurately model more realistic problems, one must be able to use mimetic operators on non-trivial regions. These mimetic operators and the \emph{mapping method} have been combined to obtain high-order accurate discretizations on non-uniform grids in \cite{Hyman, Batista}. Furthermore, in \cite{Abouali}, modified gradient and divergence operators are presented to solve Poisson's equation over a 2D curvilinear staggered grid with Robin's boundary conditions. In addition, energy conservation has been demonstrated for the 3D advection equation \cite{Dumett-Castillo-Advection3D}. For a more extensive discussion on employing high-order mimetic operators on non-trivial problems, see \cite{velazcophd2021high}.

The closest method to the mimetic operators in our investigation is the so-called summation-by-parts (SBP) operators with norms of the non-diagonal form \cite{Kreiss} which has a lower order of approximation at the boundaries than the interior of the problem domain. Moreover, \citet{Svard} showed that, in the case of SBP operators, it is impossible, to construct a coordinate transformation operator without decreasing the order of accuracy.

This paper is organized as follows: in section 2, we present a brief review of these mimetic operators, staggered grids as well as the basic mimetic operators in one dimension. The extensions to 2D and 3D can be done using Kronecker products \cite{kronecker}. These mimetic operators operators are the ones from Corbino-Castillo \cite{Corbino}, which unlike Castillo-Grone's, have no free parameters and minimal bandwidth. In section 3, it is shown that the new mimetic operators in curvilinear coordinates satisfy a discrete version of the extended Gauss divergence theorem while maintaining the same order of approximation at the boundaries and the interior of the problem domain. Section 4 presents the acoustic wave equation as one of our model problems. In section 5, we demonstrate energy conservation for the acoustic wave equation in curvilinear coordinates. In section 6, mass conservation is achieved for the acoustic wave problem in curvilinear coordinates. In section 7, numerical steady-state examples are exhibited. In section 8, we show our numerical results for our acoustic wave problem, and finally, in section 9, a summary of our conclusions and future work are presented.

\section{Mimetic Operators}

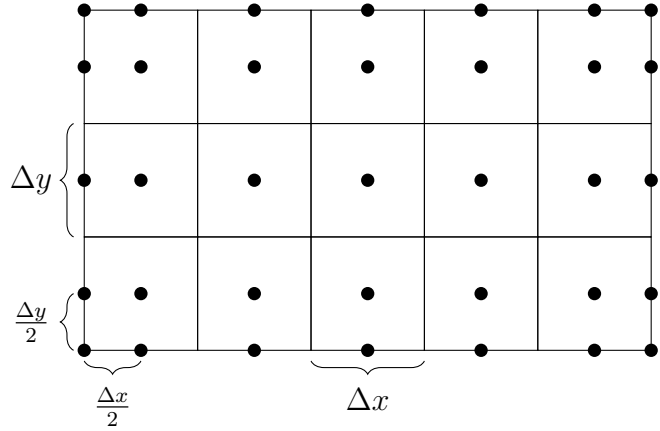
\begin{figure}[H]
		\centering
		\tikzmath{\x1 = 1.5; \r1 = 0.3;}
		\begin{tikzpicture}
			\tikzstyle{filld} = [circle, minimum width=5pt, fill, inner sep=0pt]
			\tikzstyle{diam} = [diamond, anchor=center, minimum width=0.7cm, minimum height=0.7cm, fill, inner sep=0pt]
			\draw (0,0) rectangle (\x1,\x1);
			\draw (\x1,0) rectangle (\x1*2,\x1);
			\draw (\x1*2,0) rectangle (\x1*3,\x1);
			\draw (\x1*3,0) rectangle (\x1*4,\x1);
			\draw (\x1*4,0) rectangle (\x1*5,\x1);
			\draw (0,\x1) rectangle (\x1,\x1*2);
			\draw (\x1,\x1) rectangle (\x1*2,\x1*2);
			\draw (\x1*2,\x1) rectangle (\x1*3,\x1*2);
			\draw (\x1*3,\x1) rectangle (\x1*4,\x1*2);
			\draw (\x1*4,\x1) rectangle (\x1*5,\x1*2);
			\draw (0,\x1*2) rectangle (\x1,\x1*3);
			\draw (\x1,\x1*2) rectangle (\x1*2,\x1*3);
			\draw (\x1*2,\x1*2) rectangle (\x1*3,\x1*3);
			\draw (\x1*3,\x1*2) rectangle (\x1*4,\x1*3);
			\draw (\x1*4,\x1*2) rectangle (\x1*5,\x1*3);
			\node[filld]  at (0,0) {};
			\node[filld]  at (\x1*0.5,0) {};
			\node[filld]  at (\x1*1.5,0) {};
			\node[filld]  at (\x1*2.5,0) {};
			\node[filld]  at (\x1*3.5,0) {};
			\node[filld]  at (\x1*4.5,0) {};
			\node[filld]  at (\x1*5,0) {};
			\node[filld]  at (0,\x1/2) {};
			\node[filld]  at (\x1*0.5,\x1/2) {};
			\node[filld]  at (\x1*1.5,\x1/2) {};
			\node[filld]  at (\x1*2.5,\x1/2) {};
			\node[filld]  at (\x1*3.5,\x1/2) {};
			\node[filld]  at (\x1*4.5,\x1/2) {};
			\node[filld]  at (\x1*5,\x1/2) {};
			\node[filld]  at (0,\x1*1.5) {};
			\node[filld]  at (\x1*0.5,\x1*1.5) {};
			\node[filld]  at (\x1*1.5,\x1*1.5) {};
			\node[filld]  at (\x1*2.5,\x1*1.5) {};
			\node[filld]  at (\x1*3.5,\x1*1.5) {};
			\node[filld]  at (\x1*4.5,\x1*1.5) {};
			\node[filld]  at (\x1*5,\x1*1.5) {};
			\node[filld]  at (0,\x1*2.5) {};
			\node[filld]  at (\x1*0.5,\x1*2.5) {};
			\node[filld]  at (\x1*1.5,\x1*2.5) {};
			\node[filld]  at (\x1*2.5,\x1*2.5) {};
			\node[filld]  at (\x1*3.5,\x1*2.5) {};
			\node[filld]  at (\x1*4.5,\x1*2.5) {};
			\node[filld]  at (\x1*5,\x1*2.5) {};
			\node[filld]  at (0,\x1*3) {};
			\node[filld]  at (\x1*0.5,\x1*3) {};
			\node[filld]  at (\x1*1.5,\x1*3) {};
			\node[filld]  at (\x1*2.5,\x1*3) {};
			\node[filld]  at (\x1*3.5,\x1*3) {};
			\node[filld]  at (\x1*4.5,\x1*3) {};
			\node[filld]  at (\x1*5,\x1*3) {};

			\draw [decorate, decoration={brace,mirror, amplitude=5pt, raise=4pt}] (0,0) -- (\x1*0.5,0) node [pos=0.5, below=10pt, black] {$\frac{\Delta x}{2}$};
			
			\draw [decorate, decoration={brace,mirror, amplitude=5pt, raise=4pt}] (\x1*2,0) -- (\x1*3,0) node [pos=0.5, below=10pt, black] {$\Delta x$};
			
			\draw [decorate, decoration={brace, amplitude=5pt, raise=4pt}] (0,0) -- (0,\x1*0.5) node [pos=0.5, xshift=-20pt, black] {$\frac{\Delta y}{2}$};
			
			\draw [decorate, decoration={brace, amplitude=5pt, raise=4pt}] (0,\x1) -- (0,\x1*2) node [pos=0.5, xshift=-20pt, black] {$\Delta y$};
		\end{tikzpicture}
		\caption[2D uniform staggered grid]{A 2D uniform staggered grid with $5$ horizontal cells and $3$ vertical cells. Scalar values (\tikz\draw[fill] (0,0) circle (0.3ex);), and vector components $\langle u, v \rangle$ are defined at horizontal middle cell faces and vertical middle cell faces, respectively. Observe that discrete scalar and vector fields are not defined at nodes.}
		\label{fig:staggeredGrid}
	\end{figure}

Mimetic operators, divergence ($D\equiv\nabla\cdot$), gradient ($G\equiv\nabla$), curl ($C\equiv\nabla\times$), and laplacian ($L\equiv\nabla^2$) are defined over \emph{staggered grids}. Fig.~\ref{fig:staggeredGrid} shows how a two-dimensional uniform staggered grid looks like for a resolution of $m=5$ (horizontal) by $n=3$ (vertical) cells.

Expressions (Eq.~\ref{eq:divergence4}) and (Eq.~\ref{eq:gradient4}), show the upper-left corner of the $4^{th}$-order 1D Corbino-Castillo operators, when there are $m$ cells.
\begin{align}\label{eq:divergence4}
\textit{D} = \frac{1}{\Delta x}
\begin{bmatrix}
       0        &  \cdots &       &                                                     \\\\
\frac{-11}{12}  &  \frac{17}{24}  &  \frac{3}{8}  &  \frac{-5}{24}  &  \frac{1}{24}  &    0    &  \cdots \\\\
\frac{1}{24}    &  \frac{-9}{8}   &  \frac{9}{8}  &  \frac{-1}{24}  &         0      &  \cdots           \\\\
       0        &  \ddots         &  \ddots       &  \ddots         &  \ddots        &  \ddots
\end{bmatrix}_{(m+2),(m+1)},
\end{align}

\begin{align}\label{eq:gradient4}
\textit{G} = \frac{1}{\Delta x}
\begin{bmatrix}
\frac{-352}{105} & \frac{35}{8} & \frac{-35}{24} & \frac{21}{40} & \frac{-5}{56} & 0 & \cdots \\\\
\frac{16}{105} & \frac{-31}{24} & \frac{29}{24} & \frac{-3}{40} & \frac{1}{168} & 0 & \cdots \\\\
0 & \frac{1}{24} & \frac{-9}{8} & \frac{9}{8} & \frac{-1}{24} & 0 & \cdots \\\\
\vdots & 0 & \ddots & \ddots & \ddots & \ddots
\end{bmatrix}_{(m+1),(m+2)}.
\end{align}

These discrete operators satisfy the following basic identities from vector calculus when applied to scalar fields $f, g$ and vector field $\vec v$, 
\begin{align}
  \textit{G}f_{const} &= 0, \\
  \textit{D}\vec v_{const} &= 0, \\
  \textit{C}\textit{G}f &= 0, \\
  \textit{D}\textit{C}\vec v &= 0, \\
  \textit{D}\textit{G}f &= \textit{L}f, \\
\textit{G}(fg) &= f \textit{G}g + g \textit{G}f, \\
\textit{D}(f\vec v) &= \textit{G}f \cdot \vec v + f \textit{D} \vec v, \\
\textit{L}(fg) &= f \textit{L}g + 2 \textit{G}f \cdot \textit{G}g + g \textit{L}f.
\end{align}

while maintaining a uniform order of accuracy in the entire domain. Furthermore, they satisfy the \emph{integration by parts} formula in the discrete sense,
    \begin{equation}
    h \, \langle D\vec v,f \rangle_Q + h \, \langle Gf,\vec v \rangle_P = f_n \vec v_n - f_0 \vec v_0
    \end{equation}

where $Q$ and $P$ are diagonal positive definite matrices. Higher dimensional mimetic operators are built by means of \emph{Kronecker} products \cite{Miranda}, and it can be shown that they satisfy the extended version of the Gauss divergence theorem in logical coordinates \cite{Dumett-Castillo-Identities}.

\section{Discrete Extended Gauss-Divergence Theorem in Curvilinear Coordinates}

This section shows that mimetic operators in curvilinear coordinates are conservative. Consider the discrete extended Gauss-Divergence theorem:

\begin{equation}
\left \langle D\vec v, f \right \rangle_{Q} + \left \langle \vec v, Gf \right \rangle_{P} = \left \langle B\vec v, f \right \rangle \label{eq:3}
\end{equation}

 where $B$ is known as the \emph{boundary operator}, and is defined as $B = QD + G^{T}P$. The mapping from $(x,y)$ coordinates in a physical domain $\Omega$ to $(\xi_{1}, \xi_{2})$ in a logical domain has the \emph{Jacobian} $J$ defined by:
 \begin{equation}
 J=\begin{vmatrix}
 \frac{\partial x}{\partial \xi_{1}} & \frac{\partial y}{\partial \xi_{1}} \\ 
\\
 \frac{\partial x}{\partial \xi_{2}} & \frac{\partial y}{\partial \xi_{2}} 
\end{vmatrix}
= J \left ( \xi_{1}, \xi_{2}\right )
\end{equation}


Let $\tilde{D}$ and $\tilde{G}$ denote the divergence and gradient in curvilinear coordinates, respectively. Then, the composition with the mapping coordinates results in

\begin{align}
\tilde{D} = J_{D}D,\label{eq:33} \\
\tilde{G}= J_{G}G. \label{eq:44}
\end{align}

where $J_{D}$ is the Jacobian corresponding to the divergence and $J_{G}$ is the Jacobian corresponding to the gradient.

Using \eqref{eq:33} and \eqref{eq:44} in \eqref{eq:3}, one obtains:
\begin{equation}
    \left \langle  J_{D}^{-1} \tilde{D}\vec v, f \right \rangle_{Q} + \left \langle \vec v, J_{G}^{-1} \tilde{G} f \right \rangle_{P} = \left \langle B\vec v, f \right \rangle
\end{equation}

In \cite{Batista}, it is shown that $Q^{T} = J_{D} Q_{cc}$ and $P^{T} = J_{G} P_{cc}$, where $Q_{cc}$ and $P_{cc}$ are the corresponding weights for the operators in curvilinear coordinates.
\\

Therefore,
\begin{equation}
    \left \langle  J_{D}^{-1} \tilde{D}\vec v, f  \right \rangle_{\left [J_{D}Q_{cc}\right ]^{T}} + \left \langle \vec v, J_{G}^{-1}  \tilde{G} f \right \rangle_{\left [J_{G}P_{cc}\right ]^{T}} = \left \langle B \vec v, f \right \rangle,
\end{equation}

and it follows that,
\begin{equation}
    \left \langle Q_{cc} \; J_{D}^{T} \; J_{D}^{-1} \tilde{D} \vec v,f \right \rangle + \left \langle \vec v, \; P_{cc} \;J_{G}^{T}  \; J_{G}^{-1} \tilde{G} f \right \rangle = \left \langle B \vec v, f \right \rangle,
\end{equation}

or equivalently,
\begin{equation}
    \left \langle  \tilde{D} \vec v, f\right \rangle_{Q_{cc}} + \left \langle \vec v, \tilde{G} f \right \rangle_{P_{cc}} = \left \langle B \vec v, f \right \rangle,
\end{equation}

Since $B = B_{cc}$, as was demonstrated in \cite{Batista}, then,
\begin{equation}
    \left \langle  \tilde{D} \vec v, f\right \rangle_{Q_{cc}} + \left \langle \vec v, \tilde{G} f  \right \rangle_{P_{cc}} = \left \langle B_{cc} \vec v, f\right \rangle.
\end{equation}

Ergo, we have demonstrated that the mimetic operators $\tilde{D}, \tilde{G}$ in curvilinear coordinates satisfy the extended Gauss-Divergence theorem and the corresponding \emph{exactness condition} for the curvilinear discrete divergence, using the appropriate interpolation matrices \cite{Anand-Mimetic-RRK}.

\section{Acoustic Wave Problem and Curvilinear Coordinates}
\label{S:2}

Consider a 2D spatial domain $\Omega$, so that its corresponding computational domain is described by a boundary-conforming 2D-staggered grid.

Take as a continuous model, the linearized 2D-acoustic wave motion for the scalar pressure field $p$ and the vector velocity field $\vec{V} = (u,v)$, where the scalar components of $\vec{V}$ are referred to an orthonormal unit vector basis $\left \{ \hat{e}_{1}, \hat{e}_{2} \right \}$, that is, $\vec{V}= u \hat{e}_{1} + v \hat{e}_{2}$.

The governing system of equations for a fluid having a reference unit density and constant Bulk modulus $B=(1/2)$ is:

\begin{align}
    \frac{\partial p}{\partial t} &= \left(-\frac{1}{2}\right) \nabla\cdot \vec{V}, \\
    \frac{\partial \vec{V}}{\partial t} &= -\nabla p,
\end{align}


subject to an homogeneous Dirichlet boundary condition, $p = 0$ on $\partial \Omega$, or a closed walls domain, where $ \left \langle  \vec{V}, \hat{n} \right \rangle = 0$ with $\hat{n}$ denoting the outwards unit normal to the boundary.

In \cite{boada2020high2}, it was obtained a fully discrete solution of the above linearized 3D-acoustic model, that was high-order and time stable, but using Cartesian coordinates.

There, the introduction of the adiabatic bulk modulus $B$ and a reference constant fluid density $\rho_{0}$, yielded a linearized equation of state between the pressure $p$ relative to atmospheric absolute pressure, and density $\rho$:  

\begin{equation}
\frac{\partial p}{\partial t} = \left ( \frac{B}{\rho_{0}} \right ) \; \frac{\partial \rho}{\partial t}.
\end{equation}

Analogously, the linearized \emph{continuity equation} or mass conservation became:

\begin{equation}
\rho_{0} \nabla\cdot \vec{V} =- \;  \frac{\partial \rho}{\partial t}.
\end{equation}

The particular choice $B = \frac{1}{2} \; \textup{and} \; \rho_{0} =1$, allows the use of a nice 2D - manufactured solution for testing purposes, but on the unit square as a spatial domain, covered with a uniform Cartesian staggered grid, so that the linearized continuity equation becomes:

\begin{equation}
\nabla\cdot \vec{V} = -\frac{\partial \rho}{\partial t}.   
\end{equation}

Semi-discrete energy conservation was established for the above acoustic model under Cartesian formulation in \cite{boada2020high2}.

The extension of that result, as well as mass conservation on curvilinear staggered grids that are boundary-conforming, will now be derived using the \emph{mapping method}.

\section{Semi-discrete Energy Conservation}
\label{sec:energy-conservation}

Consider the first-order energy conservation system for $p_{i}$ and $\vec{V_{i}}$, (using the standard multi-index notation for ``$i$''), and the symbol `` $\dot{}$ '' for time differentiation:

\begin{align}\label{eq: 1}
\dot{p_{i}} = \left ( - \frac{1}{2J} \right ) \tilde{D} \;\vec{V_{i}},\\ 
\dot{\vec{V_{i}}} = \left ( - \frac{1}{J} \right )  \tilde{G} \; p_{i}.
\end{align}

Here, the discrete divergence and gradient in curvilinear coordinates are defined by \eqref{eq:33} and \eqref{eq:44}.\newline
\\

From
\begin{equation}
\tilde{D}  \left ( p_{i} \vec{V_{i}} \right ) = p_{i} \; \tilde{D}  \; \left ( \vec{V_{i}} \right ) + \left \langle \tilde{G}  p_{i} , \vec{V_{i}}\right \rangle,
\end{equation}

it follows that,
\begin{equation}
    p_{i} \; \tilde{D} \left ( \vec{V_{i}} \right ) = \tilde{D} \left ( p_{i} \vec{V_{i}} \right ) - \left \langle \tilde{G}  p_{i} , \vec{V_{i}}\right \rangle = \tilde{D}  \left ( p_{i} \vec{V_{i}} \right ) + J \left \langle \dot{\vec{V_{i}}},  \vec{V_{i}} \right \rangle,
\end{equation}

and, on the other hand, 
\begin{equation}
    p_{i}  \;\dot{p_{i}} = \left ( \frac{1}{2J} \right ) \; p_{i} \; \tilde{D} \left ( \vec{V_{i}} \right ), \; \textup{or, } \; 
 p_{i} \; \tilde{D} \left ( \vec{V_{i}} \right )= \left ( -J \right ) \frac{d}{dt} \left ( p_{i}^{2} \right ),
\end{equation}

thus,
\begin{equation}
    \left ( -J \right ) \frac{d}{dt} \left ( p_{i}^{2} \right )= \tilde{D} \; \left ( p_{i} \vec{V_{i}} \right ) + J \left ( \frac{1}{2} \right ) \frac{d}{dt} \left \langle \vec{V_{i}},  \vec{V_{i}} \right \rangle,
\end{equation}

which implies,
\begin{equation}
    \frac{d}{dt} \left \{ \left ( p_i \right )^{2} + \frac{1}{2} \left \langle \vec{V_{i}},  \vec{V_{i}} \right \rangle  \right \} = - \frac{1}{J} \tilde{D}\left ( p_{i} \vec{V_{i}} \right ). \label{eq: 2}
\end{equation}

An \emph{exactness condition} for our two-dimensional operator $\tilde{D}$, has been extended from the one-dimensional case \eqref{eq: 2}, and has corresponding positive weights $M_{i}$, with $i$ multi-index.

Multiplying both sides of equation \eqref{eq: 2} by $M_{i} J \Delta x \Delta y$, and then summing over ``$i$'' to cover all cells, to get the time-derivative of the total weighted energy $E$:

\begin{equation} \label{eq1}
\begin{split}
\frac{dE}{dt} & = \frac{d}{dt} \left [ \sum_{i} M_{i} \left \{ \left ( p_{i} \right )^{2} + \frac{1}{2}  \left \langle \vec{V_{i}},  \vec{V_{i}} \right \rangle  \right \} J \left ( \Delta x \Delta y \right )\right ] \\
             & = - \left ( \Delta x \Delta y \right ) \sum_{i} M_{i} \tilde{D}\left ( p_{i} \vec{V_{i}} \right )
\end{split}
\end{equation}

The exactness condition for $\tilde{D}$ now implies $\frac{dE}{dt} = 0$ for all $t \geq  0$ when $p_{i} (t) = 0$ on $\partial\Omega$.

\section{Semi-discrete Mass Conservation}

Let us briefly review the corresponding mass conservation for the continuous linearized acoustic model:

\begin{equation}
    \frac{\partial p}{\partial t} =- p_{0} \nabla\cdot \vec{V} =-\nabla\cdot \vec{V}. 
\end{equation}

Integrating, it follows that,
\begin{equation}
\int_{\Omega} \; \frac{\partial p}{\partial t} \; dV = \frac{d}{dt} \int_{\Omega} p \; dV =- \int_{\Omega} \nabla\cdot \vec{V} \; dV 
  = \int_{\partial \Omega} \; \left \langle  \vec{V}, \hat{n} \right \rangle dS = 0,
   \label{eq:continous-model}
\end{equation} 
for a closed wall domain, where $\left \langle  \vec{V}, \hat{n} \right \rangle \; = 0 \; \textup{on} \; \partial \Omega$.

With the same notations used for the energy conservation in semi-discrete form, one would have $\dot{p}_{i} =- D \left ( \vec{V}_{i} \right )$.

After multiplying both sides of equation \eqref{eq:continous-model} by the weight $M_{i}$, assuming two dimensions, and summing over ``$i$'' to cover all cells, one gets the semi-discrete mass conservation expression:

\begin{equation}
    \frac{d}{dt} \left [ \sum M_{i} \; p_{i} \; \Delta \xi_{1} \Delta \xi_{2}\right ] = 0
\end{equation}

\section{Steady State Numerical Examples}
	\subsection{2D Poisson Equation on a Semi-Annular Region}
	
	In this example, consider the two-dimensional Poisson problem with Dirichlet boundary conditions presented in \cite{Abouali}
	
\begin{equation}
-\nabla \cdot (\nabla u) = f, \;\;\;\;(x,y) \in \Omega,
\end{equation}
where the analytical solution is,
\begin{equation}
u(x,y) = \sin{\sqrt{ x^2 + y^2 }},
\end{equation}
and the forcing term is,
\begin{equation}
f(x,y) = \frac{ \cos{\sqrt{x^2+y^2}}} { \sqrt{x^2+y^2} } - \sin{\sqrt{x^2+y^2} }.
\end{equation}

	\begin{figure}[H]
		\centering
		\begin{minipage}{0.45\textwidth}
			\centering
			\includegraphics[width=.9\linewidth]{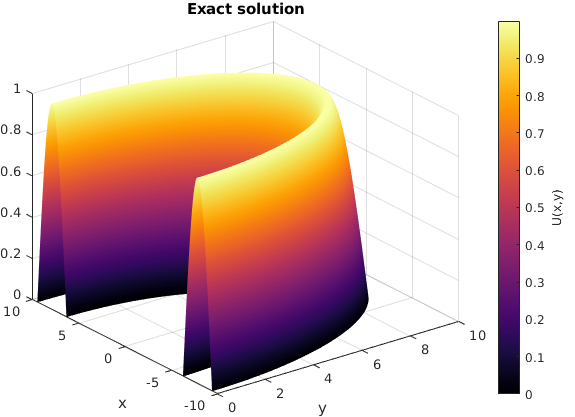}
			\caption{Exact solution for 2D Poisson equation on semi-annular region.}
			\label{fig:p2exa}
		\end{minipage}\hfill
		\begin{minipage}{0.45\textwidth}
			\centering
			\includegraphics[width=.9\linewidth]{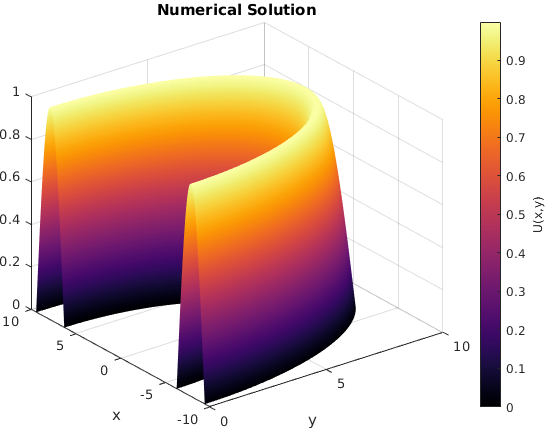}
			\caption{Computed solution for 2D Poisson equation on semi-annular region.}
			\label{fig:p2num}
		\end{minipage}
	\end{figure}
	
	\begin{table}[H]
		\centering
		\begin{tabular}{@{}lllll@{}}
			\toprule
			m   & Error($L_2$) & Order($L_2$) & Error($max$) & Order($max$) \\ \midrule
			20  & 4.9121E-06   & --           & 1.1035E-05   & --           \\
			40  & 2.6965E-07   & 4.19         & 4.6489E-07   & 4.56         \\
			80  & 1.5696E-08   & 4.10         & 2.3401E-08   & 4.31         \\
			160 & 9.4721E-10   & 4.05         & 1.4261E-09   & 4.04         \\
			320 & 5.8135E-11   & 4.03         & 8.7943E-11   & 4.02         \\ \bottomrule
		\end{tabular}
		\caption{Accuracy test for 2D Poisson equation on semi-annular region.}
		\label{tab:2dpoissoncurv}
	\end{table}
	
	The physical domain $\Omega$ is defined on a semi-annular region with polar boundaries $(r,\theta) \in [2\pi,3\pi]\times[0,\pi]$.
	
	Figures \ref{fig:p2exa} and \ref{fig:p2num} depict exact and computed solutions using $4^{th}$-order mimetic operators with a resolution of 320 cells. In Table \ref{tab:2dpoissoncurv}, it can be observed that as the resolution increases, the approximation error decreases in accordance with the anticipated convergence rate.

 \subsection{2D Poisson equation on a sinusoidal region}
	
	Examine now the following 2D Poisson problem with Dirichlet boundary conditions presented in \cite{shashkov}:
	
	\begin{equation}
	    -\nabla \cdot (\mathbf{K} \; \nabla u) = f, \;\;\;\;(x,y) \in \Omega,
	    \label{eq:2dpoisson}
	\end{equation}

	
	with anisotropic tensor $\mathbf{K}$. 
 	
	\begin{figure}[H]
		\center
		\includegraphics[scale=.5]{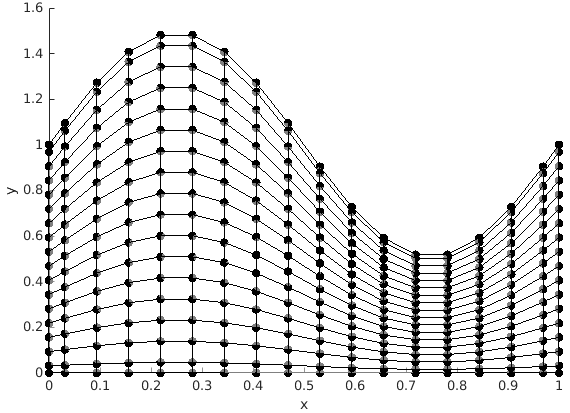} 
		\caption{Grid in domain with curvilinear boundaries.}\label{fig:stggcurvbdry}
	\end{figure}
	
	\begin{figure}[H]
		\centering
		\begin{minipage}{0.45\textwidth}
			\centering
			\includegraphics[width=.9\linewidth]{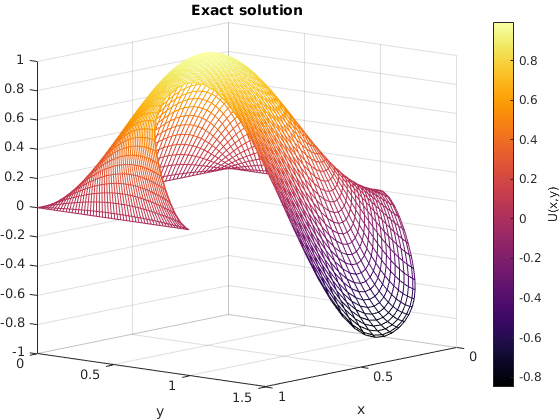}
			\caption{Exact solution for 2D Poisson equation with curvilinear boundaries.}
			\label{fig:p3exa}
		\end{minipage}\hfill
		\begin{minipage}{0.45\textwidth}
			\centering
			\includegraphics[width=.9\linewidth]{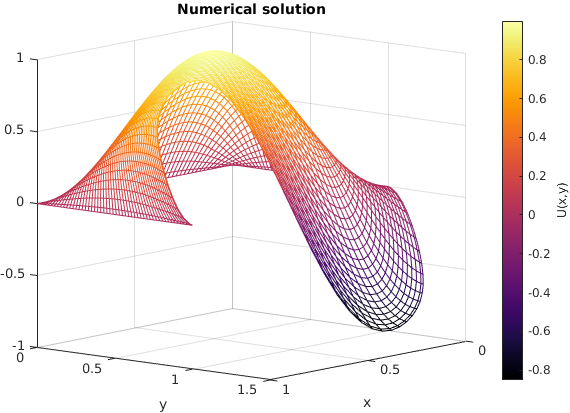}
			\caption{Computed solution for 2D Poisson equation with curvilinear boundaries.}
			\label{fig:p3num}
		\end{minipage}
	\end{figure}

 For this problem, the manufactured solution is $u = \sin(\pi x) \sin(\pi y)$, and the tensor $\mathbf{K}$ was chosen to be the result of rotating the matrix:
	\begin{equation}
	D = \begin{bmatrix}
        d_1 & 0\\ 
        0 & d_2
        \end{bmatrix},
	\end{equation}
    therefore, $\mathbf{K} = RDR^T$, where:
	\begin{equation}
	    R = \begin{bmatrix}
        \cos(\theta) & -\sin(\theta)\\ 
        \sin(\theta) & \;\;\;\cos(\theta)
    \end{bmatrix},
    \end{equation}
	and,
 \begin{equation}
     \theta = \pi / 4.
 \end{equation}
 
	The spatial discretization for this problem is shown in figure \ref{fig:stggcurvbdry}. Finally, Figures \ref{fig:p3exa} and \ref{fig:p3num} depict exact and computed solutions using $4^{th}$-order mimetic operators with a resolution of $25 \times 25$ cells. 

Table \ref{tab:2dpoissoncurv2} contains the accuracy test for this problem. Results in this table show that the $4^{th}$-order curvilinear operators attain the targeted order of accuracy.

	\begin{table}[H]
		\centering
		\begin{tabular}{@{}lllll@{}}
			\toprule
			m   & Error($L_2$) & Order($L_2$) & Error($max$) & Order($max$) \\ \midrule
			20  & 1.1987E-05   & --            & 4.3227E-05   & --      \\
			40  & 8.3285E-07   & 3.8473        & 3.0574E-06   & 3.8216  \\
			80  & 5.4609E-08   & 3.9308        & 1.9955E-07   & 3.9375  \\
			160 & 3.4870E-09   & 3.9691        & 1.2660E-08   & 3.9785  \\
			320 & 2.2011E-10   & 3.9857        & 7.9526E-10   & 3.9927  \\ \bottomrule
		\end{tabular}
		\caption{2D Poisson equation on a sinusoidal region.}
		\label{tab:2dpoissoncurv2}
	\end{table}
	
	\subsection{3D Poisson Equation in a Sinusoidal Volume}
	
	In this case, the Poisson problem in a 3D sinusoidal volume is solved,
	
	\begin{equation}
	    \nabla^2 u = f, \;\;\;\;(x,y) \in \Omega
	    \label{eq:3dpoisson}
	\end{equation}
	
	with exact solution $u = x^2+y^2+z^2$, and Dirichlet boundary conditions.
	
\begin{figure}[H]
		\centering
		\begin{minipage}{0.45\textwidth}
			\centering
			\includegraphics[width=.9\linewidth]{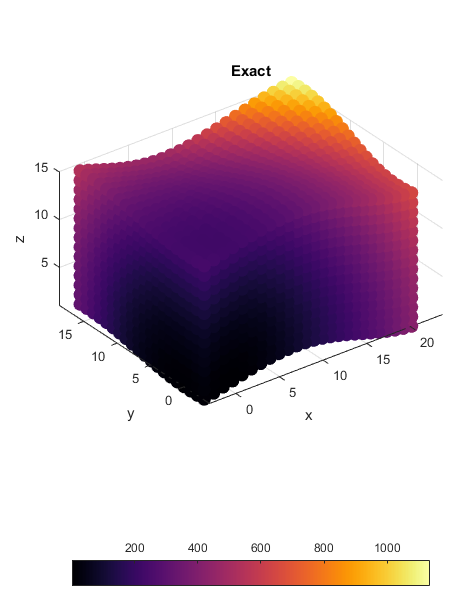}
			\caption{Exact solution for the 3D Poisson problem.}
			\label{fig:3dexact}
		\end{minipage}\hfill
		\begin{minipage}{0.45\textwidth}
			\centering
			\includegraphics[width=.9\linewidth]{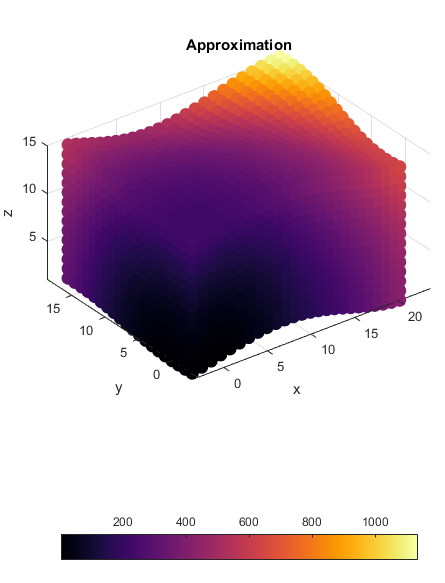}
			\caption{Computed solution for the 3D Poisson problem.}
			\label{fig:3dapprox}
		\end{minipage}
	\end{figure}
Figures \ref{fig:3dexact} and \ref{fig:3dapprox} illustrate the exact and computed solution for this problem using the $4^{th}$-order 3D curvilinear operators. In addition, its corresponding accuracy test can be found in Table \ref{tab:3dpoissoncurv}. As shown in the previous 2D problems, the desired accuracy order is also attained in the 3D case.

\begin{table}[H]
		\centering
		\begin{tabular}{@{}lllll@{}}
			\toprule
			m   & Error($L_2$) & Order($L_2$) & Error($max$) & Order($max$) \\ \midrule
			20  & 4.9121E-06   & --           & 1.1035E-05   & --           \\
			40  & 2.6965E-07   & 4.19         & 4.6489E-07   & 4.56         \\
			80  & 1.5696E-08   & 4.10         & 2.3401E-08   & 4.31         \\
			160 & 9.4721E-10   & 4.05         & 1.4261E-09   & 4.04         \\
			320 & 5.8135E-11   & 4.03         & 8.7943E-11   & 4.02         \\ \bottomrule
		\end{tabular}
		\caption{3D Poisson equation on a sinusoidal volume.}
		\label{tab:3dpoissoncurv}
	\end{table}
 
\section{Acoustic Wave Numerical Example}
\begin{figure}[H]
		\center
		\includegraphics[scale=.5]{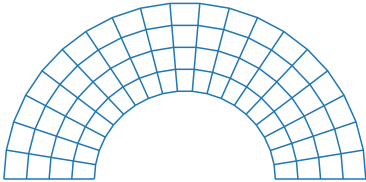} 
		\caption{Region $\Omega$: a $\leq$ r $\leq$ b  semi-annulus.}\label{fig:semi-annulus}
	\end{figure}

	In this problem, the focus is on the 2D acoustic wave equation on a semi-annular region $\Omega$ depicted in figure \ref{fig:semi-annulus}. 
 
	\begin{figure}[H] 
		\begin{subfigure}[b]{0.5\linewidth}
			\centering
			\includegraphics[width=0.9\linewidth]{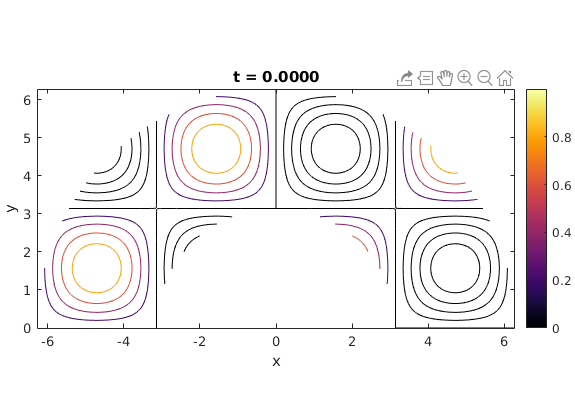} 
			\caption{Initial condition} 
			\label{fig7:a} 
		\end{subfigure}
		\begin{subfigure}[b]{0.5\linewidth}
			\centering
			\includegraphics[width=0.9\linewidth]{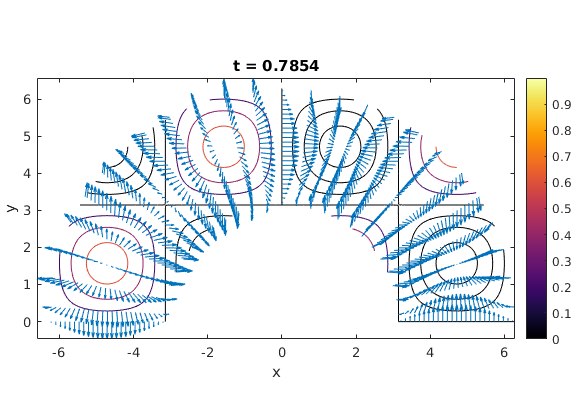} 
			\caption{Numerical solution at $t = \pi/4$} 
			\label{fig7:b} 
		\end{subfigure} 
		\begin{subfigure}[b]{0.5\linewidth}
			\centering
			\includegraphics[width=0.9\linewidth]{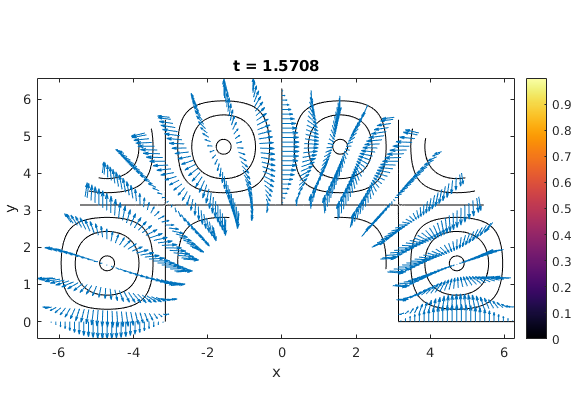} 
			\caption{Numerical solution at $t = \pi/2$} 
			\label{fig7:c} 
		\end{subfigure}
		\begin{subfigure}[b]{0.5\linewidth}
			\centering
			\includegraphics[width=0.9\linewidth]{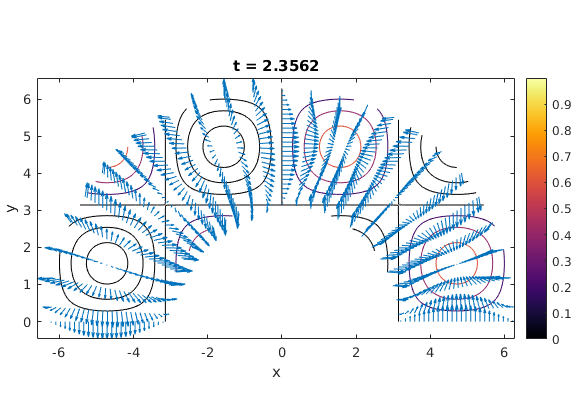} 
			\caption{Numerical solution at $t = 3\pi/2$} 
			\label{fig7:d} 
		\end{subfigure} 
		\caption{Pressure contour plot, and vector field (velocity) visualization at four different times.}
		\label{fig:linwave1} 
	\end{figure}
	
Consider that $p = p(x,y) \; \textup{and} \; \vec{V}= (u,v) \;\textup{satisfy}$:
	
	\begin{equation}
	    \frac{\partial p}{\partial t} = \left ( - \frac{1}{2} \right ) \nabla\cdot \vec{V},
     \end{equation}
	
	\begin{equation}
	\frac{\partial \vec{V}}{\partial t} =- \nabla p.
    \end{equation}

	for these, the exact solutions are:
	
	\begin{align}
	p &= \sin  x \; \sin  y \; \cos   t,\\ 
	u &= -\cos  x \; \sin  y \; \sin t,\\ 
	v &= - \sin  x \; \cos  y \; \sin  t.
	\end{align}
	
	Each scalar component of $\vec{V}$ satisfies the scalar wave equation. However,
	
	$\frac{\partial u}{\partial x}= \; \sin  x \; \sin y \; \sin t \; \equiv 0$, while $u = - \; \cos  x \; \sin  y \; \sin t \; \equiv 0$, when $\sin  y = 0$, so that it is a Robin boundary condition.
	
	Figure \ref{fig:linwave1} depicts contour plots and the corresponding vector field visualization for the numerical solution computed using mimetic operators at different times. 
 
 Table \ref{tab:linwave} displays the numerical error in $L_2$-norm for the computed solution of $p$, $u$, and $v$. Here, the high-order convergence demonstrates the mimetic approach's efficiency in a non-rectangular geometry.
		
	\begin{table}[H]
		\centering
		\begin{tabular}{@{}lll|ll|ll@{}}
			\toprule
			\multirow{2}{*}{m} & \multicolumn{2}{c}{$err(p)$} & \multicolumn{2}{c}{$err(u)$} & \multicolumn{2}{c}{$err(v)$} \\ \cmidrule(l){2-7} 
			& Error($L_2$)  & Order($L_2$) & Error($L_2$)  & Order($L_2$) & Error($L_2$)  & Order($L_2$) \\ \midrule
			16                 & 3.7483E-03    & --           & 2.4710E-03    & --           & 2.1802E-03    & --         \\
			32                 & 3.4178E-04    & 3.46         & 1.2705E-04    & 4.28         & 1.2136E-04    & 4.17       \\
			64                 & 2.3358E-05    & 3.87         & 1.2072E-05    & 3.40         & 7.3432E-06    & 4.05       \\
			128                & 2.2351E-06    & 3.39         & 1.2109E-06    & 3.32         & 5.1292E-07    & 3.84       \\
			256                & 2.4048E-07    & 3.21         & 1.3654E-07    & 3.15         & 3.9980E-08    & 3.68       \\ \bottomrule
		\end{tabular}
		\caption{Numerical results for $p$, $u$, $v$ at $t = 1$ using $\Delta t = 0.001$.}
		\label{tab:linwave}
	\end{table}
	
\section{Concluding Remarks}
In this study, we systematically reviewed the relevant literature on mimetic methods, setting the stage for a comprehensive exploration of their application in curvilinear staggered grids. We presented the methodology for constructing and applying mimetic operators on curvilinear staggered grids. We extended the Corbino-Castillo operators so they can be employed to solve 2D and 3D problems in non-trivial geometries. We demonstrated that our resulting curvilinear operators satisfy the extended Gauss divergence theorem in its discrete form, and we assessed for energy and mass conservation, in curvilinear coordinates, for the acoustic wave case. We conducted tests using these novel operators in diverse scenarios, and the outcomes aligned with our anticipated results. We are already working on using our mimetic operators on overlapping grids and will look at embedded boundaries afterward.
\section*{Acknowledgements}
The authors would like to thank Dr. Guillermo Miranda of Universidad Central de Venezuela for the helpful discussions on the subject covered in this paper.

\bibliographystyle{elsarticle-num-names}
\bibliography{sample.bib}

\end{document}